\documentclass[a4paper,11pt,reqno]{amsart}
\usepackage[a4paper,hscale=0.7,vscale=0.7,centering]{geometry}
\usepackage{amssymb}
\usepackage{color}

\newcommand{\E}{\mathbb{E}}
\renewcommand{\P}{\mathbb{P}}
\newcommand{\enne}{\mathbb{N}}
\newcommand{\erre}{\mathbb{R}}

\newcommand{\dom}{\mathop{\mathrm{dom}}\nolimits}
\newcommand{\eps}[1]{#1^{(\varepsilon)}}
\newcommand{\ip}[2]{\langle #1,#2 \rangle}
\newcommand{\bip}[2]{\big\langle #1,#2 \big\rangle}

\newtheorem{prop}{Proposition}
\newtheorem{thm}[prop]{Theorem}
\newtheorem{coroll}[prop]{Corollary}
\newtheorem{lemma}[prop]{Lemma}
\newtheorem{defi}[prop]{Definition}
\theoremstyle{remark}
\newtheorem{rmk}[prop]{Remark}

\newcommand{\cf}{\mathcal F}
\newcommand{\ch}{\mathcal H}

\newcommand{\cl}{\mathcal L}

\newcommand{\al}{\alpha}

\newcommand{\ep}{\varepsilon}

\begin{document}
\title[Stochastic nonlinear wave equations]{Existence of weak
  solutions for a class of semilinear stochastic wave equations}

\author{Carlo Marinelli}
\address[C.~Marinelli]{Facolt\`a di Economia, Universit\`a di Bolzano,
  I-39100 Bolzano, Italy and Dipartimento di Matematica, Universit\`a
  di Trento, I-38123 Trento, Italy.}
\urladdr{http://www.uni-bonn.de/$\sim$cm788}

\author{Llu\'is Quer-Sardanyons}
\address[L.~Quer-Sardanyons]{Departament de Matem\`atiques,
  Universitat Aut\`onoma de Barcelona, 08193 Cerdanyola del Vall\`es
  (Barcelona), Catalunya, Spain.}
\urladdr{http://www.mat.uab.cat/$\sim$quer}

\date{21 October 2011}

\begin{abstract}
  We prove existence of weak solutions (in the probabilistic sense)
  for a general class of stochastic semilinear wave equations on
  bounded domains of $\erre^d$ driven by a possibly discontinuous
  square integrable martingale.
\end{abstract}



\thanks{Part of the work for this paper was carried out while the
  first author was visiting the Department of Statistics of Purdue
  University supported by a MOIF fellowship, and while the second
  author was visiting the CIRM at the University of Trento supported
  by a fellowship of the Fondazione Bruno Kessler. The second author
  was also supported by the grant MEC-FEDER Ref.~MTM2009-08869 of the
  Ministerio de Ciencia e Innovaci\'on, Spain. We thank Viorel Barbu
  for helpful conversations on the topic of this paper.}

\maketitle

\section{Introduction}
The purpose of this paper is to prove existence of weak solutions, in
the probabilistic sense, for a class of semilinear stochastic wave
equations of the type
\begin{equation}     \label{eq:fw}
\frac{\partial^2 u}{\partial t^2}(t,x) 
- \Delta u(t,x) + \beta(u(t,x)) = \dot{\eta}(t,x)
\end{equation} 
on a bounded space-time domain $[0,T] \times D \subset \erre^{1+d}$,
complemented with suitable boundary conditions. Here $\beta$ denotes a
maximal monotone graph in $\erre \times \erre$ and $\dot{\eta}$ stands
for the time derivative of a stochastic integral with respect to a
Hilbert space-valued (possibly discontinuous) martingale which may
depend on $u$. Of course (\ref{eq:fw}) is only a formal expression,
whose corresponding rigorous stochastic evolution equation, as well as
the definition of solution, will be given below.

Existence of solutions for stochastic PDEs such as (\ref{eq:fw})
cannot be obtained (to the best of our knowledge) simply as an
application of some general technique. In particular, on the one hand
(\ref{eq:fw}) cannot be cast in the variational setting of Pardoux
\cite{Pard} and Krylov-Rozovskii \cite{KR-spde}, as the equation is
not of monotone type, in spite of $\beta$ being monotone. On the other
hand, since the nonlinear term $\beta$ does not satisfty any (even
local) Lipschitz condition, the semigroup approach does not seem to be
applicable directly either. However, major efforts have been devoted,
especially in recent years, to obtain existence, uniqueness, and
regularity results for classes of stochastic PDEs that do not fall
into any standard framework (see e.g. \cite{BDPR-porous}), and our
work is a contribution in this direction.
In the particular case of stochastic wave equations of the form
(\ref{eq:fw}), some well-posedness results in the mild sense have been
obtained assuming that $\beta$ has sufficiently slow polynomial growth
and that $\eta$ is a stochastic integral with respect to a Wiener
process  -- see e.g. \cite{Chow-wave, MM-wave, ondrejat-wave}.
Moreover, Ondrej\'at \cite{Ond:mgwave} has recently proved that for
any $\beta$ continuous and polynomially growing, and $\eta$ a
stochastic integral with respect to a spatially homogeneous Wiener
process with finite spectral measure, (\ref{eq:fw}) admits a (suitably
defined) global weak solution. More precisely, he proves the existence
of a solution to an infinite dimensional martingale problem associated
to (\ref{eq:fw}).

Let us also recall that there exist in the literature several
well-posedness results for stochastic wave equations with Lipschitz
continuous drift, among which \cite{kar-zabczyk, peszat,
  peszat-zabczyk2, peszat-zabczyk}, where the semigroup approach is
used, and \cite{carmona-nualart, conus-dalang, dalang-mini, dalang-frangos,
  dalang-mueller, millet-sanz}, that adopt the random field approach
(\`a la Walsh \cite{walsh}).

Our approach, whose origins can be traced back to \cite[\S
1.6]{Sko-studies} (cf. also \cite[Section 2.6]{krylov} and
\cite{Pell-mp}), is entirely different and relies instead on
Skorohod's representation theorem and weak convergence results for
stochastic integrals with respect to general martingales. This route
allows us to avoid going through the martingale problem for an
equation with discontinuous noise, which is already quite involved in
finite dimensions (see e.g. \cite{Jac-LNM, LepMar}). Roughly speaking,
our proof proceeds as follows: we obtain a priori estimates for
solutions of regularized equations (i.e. with smoother $\beta$) which
imply tightness in an appropriate topology. We can thus construct a
sequence of processes converging almost surely on a different
probability space. The final step consists in showing that the limit
process is a weak solution of the equation of interest, in a sense
made precise below. A key ingredient in this step is played by a
convergence result in \cite{Bre-mm} (see Theorem \ref{thm:brezis}
below). Let us also mention that unfortunately we cannot prove
well-posedness, but only existence of a solution. Nonetheless, this is
consistent with the deterministic case (i.e. with $\eta \equiv 0$),
for which, under the present generality of $\beta$, no well-posedness
results are known, even though existence of global weak solutions is
known in some cases (see e.g. \cite{Bre-mm, ShaStr-book}).

The paper is organized as follows: in Section \ref{main} we introduce
notation, recall some basic preliminaries needed throughout the paper,
define the concept of solution to (\ref{eq:fw}), and state the main
result of the paper (Theorem \ref{thm:main}). Section~\ref{sec:aux}
contains a few auxiliary results on mild solutions of stochastic
evolution equations with Lipschitz nonlinearities that are needed in
the proof of the main result and that might be of independent
interest. Finally, Section~\ref{proof:main} is devoted to the proof of
the main result.

\section{Main result}     \label{main}
\subsection{Notation and preliminaries} \label{not-prel}
Let $D \subset \erre^d$ be a bounded domain with smooth boundary
$\partial D$ and $T$ a fixed positive real number. We shall use
standard notation for spaces of integrable functions and Sobolev
spaces on $D$. In particular, $H^1_0(D)$ will denote the closure of
$C^\infty_c(D)$ in the topology of $H^1(D)$. When no confusion may
arise, we shall suppress the indication of the domain $D$, so that
$L^2$ and $H^1_0$ stand for $L^2(D)$ and $H^1_0(D)$, respectively. We
shall denote the Laplace operator on $D$ with Dirichlet boundary
conditions by $\Delta$.

Let $\beta$ be a maximal monotone graph in $\erre \times \erre$ such
that $\dom(\beta)=\erre$ and $0\in \beta(0)$. Recall that a graph
$\beta$ is called monotone if, for any $x_1$, $x_2 \in \erre$, $b_1
\in \beta(x_1)$, $b_2 \in \beta(x_2)$, one has $(b_1-b_2)(x_1-x_2)
\geq 0$, and a monotone graph $\beta$ is maximal if it is not properly
contained in any other monotone graph. Let us recall that, if
$f_0:\erre \to \erre$ is an increasing function, the graph $f: x
\mapsto [f_0(x-),f_0(x+)] \cap \erre$ is a monotone graph in $\erre
\times \erre$. As a matter of fact, all maximal monotone graphs in
$\erre \times \erre$ are constructed in this manner (cf. \cite[Exemple
2.3.1]{Bmax}).
We shall also assume, without loss of generality (see
e.g. \cite[p.~71]{barbu}), that $\beta=\partial j$ for a continuous
convex function $j:\erre\to\erre$, $j(x)\geq 0$ for all $x \in
\erre$. Here $\partial$ stands for the subdifferential in the sense of
convex analysis, i.e.
\[
\partial j(x) := \big\{
y\in \erre:\; j(z)-j(x) \geq y(z-x) \quad \forall z \in \erre
\big\}.
\]
If $j$ is differentiable at $x$, then $\partial j(x)$ reduces to a
single point and coincides with $f'(x)$. For all notions of convex
analysis and the theory of monotone operators used in the paper we
refer to \cite{barbu, Bmax}.

For any real Hilbert spaces $E$ and $F$, let us denote the space of
linear bounded operators from $E$ to $F$ by $\mathcal{L}(E \to F)$,
and its subspace of Hilbert-Schmidt operators from $E$ to $F$ by
$\mathcal{L}_2(E \to F)$. Given a symmetric nonnegative nuclear
operator  $R \in \mathcal{L}(E \to E)$, we shall use the notation
$\mathcal{L}_2^R(E \to F):=\mathcal{L}_2(R^{1/2}E \to F)$, and
$|B|_R:=|BR^{1/2}|_{\mathcal{L}_2(E \to F)}$.

Given a (fixed) real separable Hilbert space $\ch$, let $M$ be a
$\ch$-valued square integrable martingale defined on a complete
filtered probability space $(\Omega, \mathcal{F},
\{\mathcal{F}_t\}_{t\in[0,T]}, \P)$, satisfying the ``usual''
conditions. Let $Q_M$ denote the martingale covariance of $M$,
i.e. the (unique) predictable process with values in
$\mathcal{L}_1^+(\ch)$, the space of symmetric nonnegative nuclear
operators on $\ch$, such that $\langle\!\langle M
\rangle\!\rangle(t)=\int_0^t Q_M(s)d\langle M \rangle(s)$, where
$\langle M \rangle$ and $\langle\!\langle M \rangle\!\rangle$ stand
for the Meyer process and the operator angle bracket of $M$,
respectively. Here and in the following we use standard notation and
terminology for infinite dimensional stochastic calculus, for which we
refer to \cite{Met} (see also \cite{PZ-libro}).
Unless otherwise stated, we shall always assume that there exists $Q
\in \mathcal{L}_1^+(\mathcal{H})$ such that
\[
\langle\!\langle M \rangle\!\rangle (t) - \langle\!\langle M
\rangle\!\rangle (s) \leq (t-s) Q
\]
for all $0\leq s\leq t \leq T$. Recall that, if $X$ is a further
Hilbert space, then any predictable process $\Phi:[0,T] \to
\mathcal{L}_2^Q(\mathcal{H} \to X)$ is integrable with respect to $M$,
and it holds
\begin{equation}     \label{eq:16}
  \E \Big| \int_0^T \Phi(s)\,dM(s) \Big|^2_X =
  \E \int_0^T |\Phi(s)|_{Q_M}^2 \,d\langle M\rangle(s) \leq
  \E \int_0^t |\Phi(s)|^2_Q \, ds.
\end{equation}
Note that any martingale L\'evy process (in particular, a Wiener
process) with nuclear covariance operator satisfies the hypotheses
imposed on $M$.

The space of $E$-valued random variables with finite $p$-th moment
will be denoted by $\mathbb{L}^p(E)$, without explicit mention of the
underlying probability space if no confusion may arise.  Finally, the
set of c\`adl\`ag functions defined on $[0,T]$ and taking values in
$E$ (here $E$ can be any Polish space) will be denoted by $D([0,T]\to
E)$.

We shall write $a \lesssim b$ to mean that there exists a constant
$N>0$ such that $a \leq N b$.

\subsection{Main result}
Let us begin specifying the definition of solution for (\ref{eq:fw}),
which must be interpreted as the system of equations
\begin{equation}     \label{eq:sw}
\left\{
\begin{array}{l}
du(t)=v(t)\,dt,\\
dv(t)-\Delta u(t)\,dt + \beta(u(t))\,dt = G_0(u(t-))\,dM(t),
\end{array}
\right.
\end{equation}
with initial conditions $(u(0),v(0))=(u_0,v_0) \in \mathbb{L}^2(H^1_0
\times L^2)$, so that, on a formal level,
\[
\eta(t,\cdot) = \int_0^t G_0(u(s-))\,dM(s).
\]
We shall assume throughout the paper that $G_0: L^2 \to \cl_2^Q(\ch
\to L^2)$ is Lipschitz continuous.
\begin{defi}     \label{def:wsol}
  A \emph{weak solution} of equation (\ref{eq:sw}) consists of a
  stochastic basis $\bar{\mathcal{B}}:=(\bar{\Omega}, \bar{\cf},
  \bar{\mathbf{F}}, \bar{\mathbb{P}})$,
  $\bar{\mathbf{F}}:=\{\bar{\mathcal{F}}_t\}_{t\leq T}$, an $\ch$-valued
  $\bar{\mathbf{F}}$-martingale $\bar{M}$,
  and an $\bar{\mathbf{F}}$-adapted c\`adl\`ag process $\bar{u}: [0,T] \to
  L^2$, both defined on $\bar{\mathcal{B}}$, such that the
  following conditions are satisfied:
  \begin{itemize}
  \item[(i)] $\bar{M}$ has the same law of $M$ in $D([0,T] \to L^2)$;
  \item[(ii)] there exists $g\in L^1(\bar{\Omega}\times [0,T]\times
    D)$ such that $g(\bar{\omega},t,x)\in
    \beta(\bar{u}(\bar{\omega},t,x))$ for
    $\bar{\P}\otimes\mathrm{Leb}$-a.a. $(\bar{\omega},t,x)$ and
    \begin{multline}
      \bar{u}(t) + 
      \int_0^t \frac{1}{\sqrt{-\Delta}}\sin\big((t-s)\sqrt{-\Delta}\big)
               g(s)\,ds = \cos\big(t\sqrt{-\Delta}\big) u_0
      + \frac{1}{\sqrt{-\Delta}}\sin\big(t\sqrt{-\Delta}\big)v_0 \\
      + \int_0^t \frac{1}{\sqrt{-\Delta}}\sin\big((t-s)\sqrt{-\Delta}\big)
                 G_0(\bar{u}(s-))\,d\bar{M}(s)
      \label{eq:duh}
    \end{multline}
    $\bar{\P}$-a.s. for all $t \leq T$.
  \end{itemize}
\end{defi}

\begin{rmk}
  Equation (\ref{eq:duh}) is motivated by the classical Duhamel's
  representation of solutions to the linear non-homogeneous wave
  equation (see e.g. \cite[\S 4.1]{ShaStr-book}). For the right
  functional spaces in which the integrals have to be understood we
  refer to the proof of Theorem~\ref{thm:main} below.
\end{rmk}

Our main result is the following theorem, which establishes existence
of a weak solution for (\ref{eq:sw}).
\begin{thm}    \label{thm:main} 
  Assume that $u_0 \in \mathbb{L}^2(H_0^1)$, $v_0 \in
  \mathbb{L}^2(L^2)$ and $\E|j(u_0)|_{L^1}<\infty$. Then equation
  (\ref{eq:sw}) admits a weak solution.
\end{thm}


\section{Auxiliary results}     \label{sec:aux}
We collect in this section some auxiliary results on regularization
and a priori estimates for mild solutions of stochastic evolution
equations (with Lipschitz nonlinearities) which will be used in the
next section. These results are not tied in any ways to the
specific wave equation introduced above.

In addition to the notation already introduced, throughout this
section we shall fix a real separable Hilbert space $X$ and a linear
maximal monotone operator $A: \dom(A) \subset X \to X$, and we shall
set, for the sake of compactness of notation,
$\mathcal{L}_2^Q:=\mathcal{L}_2^Q(\mathcal{H} \to X)$. The norm and
the scalar product in $X$ will be denoted by $|\cdot|$ and
$\ip{\cdot}{\cdot}$, respectively. There will be no risk of confusion
with the notation used in the other sections.

The following simple result can be proved essentially as
\cite[Lemma~2.4.1]{PreRoeck}, thus we omit its proof (cf. also
\cite[p.~114]{PZ-libro}).
\begin{lemma}     \label{l3}
  Let $\Phi:[0,T] \to \cl_2^Q$ be a predictable process such that
  \begin{equation}     \label{30}
  \E \int_0^T |\Phi(s)|_Q^2 \,ds < \infty,
  \end{equation}
  and $B:X \to X$ a bounded linear operator. Then
  \[
  B \int_0^t \Phi(s)\, dM(s) = \int_0^t B \Phi(s)\, dM(s)
  \]
  $\P$-a.s. for all $t \in [0,T]$.
\end{lemma}
As an immediate consequence we also have the following corollary,
which exploits the well-known fact that $-A$ and its associated
semigroup and resolvent commute.
\begin{coroll}     \label{cor}
  Let $\Phi$ be as in Lemma \ref{l3} and set $\eps{\Phi}:= (I+\ep
  A)^{-1}\Phi$, $\varepsilon>0$. Then
  \begin{equation}    \label{15}
    A \int_0^t e^{-(t-s)A} \eps{\Phi}(s) \,dM(s) =
    \int_0^t A e^{-(t-s)A} \eps{\Phi}(s) \,dM(s)
  \end{equation}
 $\P$-a.s. for all $t \in [0,T]$.
\end{coroll}
In order to obtain a priori estimates for mild solutions we shall
approximate by strong solutions. The next lemma is a tool for the
proposition to follow. Since we have not been able to find the proofs of these results in the literature,
we will include them for the reader's convenience and the sake of completeness.
\begin{lemma}     \label{l2}
  For $\varepsilon>0$, let $\eps{Y}$ be the unique mild solution of
  the equation
  \begin{equation}
    d\eps{Y}(t) + A\eps{Y}(t)\,dt + f(t)\,dt =
    \eps{\Phi}(t)\,dM(t), \qquad \eps{Y}(0)=Y_0,
    \label{13}
  \end{equation}
  where $Y_0\in \dom(A)$, $f \in L^1([0,T] \to \dom(A))$ $\P$-a.s., and
  $\Phi$, $\eps{\Phi}$ are as in Corollary \ref{cor}.  Then $\eps{Y}$
  is a strong solution of (\ref{13}).
\end{lemma}
\begin{proof}
  As a first step, we shall assume that $f$ and $Y_0$ are identically
  zero, so that
  \begin{align}     \label{eq:14-1}
  \eps{Y}(t) & = \int_0^t e^{-(t-s)A} \eps{\Phi}(s)\,dM(s) \nonumber \\
		& =(I+\varepsilon A)^{-1} \int_0^t e^{-(t-s)A} \Phi(s)\,dM(s)
  \qquad \forall t \in [0,T].
  \end{align}
  We want to prove that one has
  \begin{equation}    \label{14}
    \eps{Y}(t)=-\int_0^t A \eps{Y}(s)\,ds + \int_0^t \eps{\Phi}(s)\,dM(s)
    \qquad \forall t\in [0,T].
  \end{equation} 
  On the one hand, it is clear that the stochastic integral on the
  right-hand side of (\ref{14}) is well-defined. On the other hand,
  since $\eps{Y}(s) \in \dom(A)$ $\P$-a.s. for all $s\in [0,t]$,
  applying a stochastic Fubini's theorem (e.g. as formulated in
  \cite{Leon-Fub}) and Corollary \ref{cor}, we have
  \begin{align*}
    \int_0^t A \eps{Y}(s) ds &= \int_0^t \left( A \int_0^s e^{-(s-r)A} \eps{\Phi}(r) dM(r)\right) ds\\
    &= \int_0^t \left( \int_0^s  A e^{-(s-r)A} \eps{\Phi}(r) dM(r)\right) ds\\
    &= \int_0^t \left( \int_r^t  A e^{-(s-r)A} \eps{\Phi}(r) ds \right) dM(r) \\
    &= \int_0^t \left( - e^{-(t-r)A} \eps{\Phi}(r) + \eps{\Phi}(r) \right) dM(r) \\
    &= - \int_0^t e^{-(t-r)A} \eps{\Phi}(r)\,dM(r) + \int_0^t \eps{\Phi}(r)\,dM(r),
  \end{align*}
  which is equivalent to (\ref{14}), in view of (\ref{eq:14-1}).

  In the general case (i.e. without assuming that $f$ and $Y_0$ are
  identically zero), the mild solution of (\ref{13}) is given by
  \[
  \eps{Y}(t) + \int_0^t e^{-(t-s)A} f(s)\,ds = e^{-tA} Y_0 
  + \int_0^t e^{-(t-s)A} \eps{\Phi}(s)\,dM(s)
  \]
  for all $t \in [0,T]$. Setting
  \begin{equation}     \label{eq:z}
  \eps{Z}(t):= \eps{Y}(t) + \int_0^t e^{-(t-s)A} f(s)\, ds - e^{-tA} Y_0,
  \end{equation}
  we have $\eps{Z}(t)\in \dom(A)$ for all $t\in [0,T]$ and
  \[
  \eps{Z}(t)=\int_0^t e^{-(t-s)A} \eps{\Phi}(s) \,dM(s)
  \]
  for all $t\in [0,T]$.
  Therefore, by the first part of the proof, $\eps{Z}$ also verifies
  \[
  \eps{Z}(t)= - \int_0^t A \eps{Z}(s) \,
  ds + \int_0^t \eps{\Phi}(r)\, dM(r) \qquad \forall t\in [0,T],
  \]
  which implies, recalling (\ref{eq:z}) and applying Fubini's theorem,
  \[
  \eps{Y}(t) + \int_0^t A \eps{Y}(s)\, ds + \int_0^t f(s)\, ds = Y_0 +
  \int_0^t \eps{\Phi}(r) dM(r),
  \]
  thus concluding the proof.
\end{proof}
The next proposition establishes a priori estimates for the mild
solution of an equation with Lipschitz nonlinearities.
\begin{prop}\label{prop3}
  Let $y$ be the mild solution of
  \[
  dy(t) + Ay(t)\,dt = f(y(t))\,dt + \Phi_0(y(t-))\,dM(t),
  \qquad y(0)=y_0,
  \]
  where $f: X \to X$ and $\Phi_0: X \to
  \cl_2^Q$ are Lipschitz continuous. Then
\begin{equation}     \label{12}
  \E \sup_{t \leq T} |y(t)|^2
  \lesssim 1 + \E |y_0|^2 + \E\sup_{t \leq T} \int_0^t \ip{f(y(s))}{y(s)} \,ds 
  + \int_0^T \big(\E \sup_{r\leq s} |y(r)|^2 \big) ds.
\end{equation}
\end{prop}
\begin{proof}
  Let us recall that the mild solution $y$ is unique, c\`adl\`ag, and
  satisfies (see e.g. \cite{Kote-Doob})
  \[
  \E \sup_{t \leq T} |y(t)|^2 < +\infty.
  \]
  Set $b(\cdot):=f(y(\cdot))$, $\Phi(\cdot):=\Phi_0(y(\cdot-))$. Then
  $y$ is the mild solution of
  \[
  dy(t) + Ay(t)\,dt = b(t)\,dt + \Phi(t)\,dM(t),
  \qquad y(0)=y_0,
  \]
  that is
  \[
  y(t) = e^{-tA}y_0 + \int_0^t e^{-(t-s)A} b(s)\,ds
  + \int_0^t e^{-(t-s)A} \Phi(s)\,dM(s)
  \]
  $\P$-a.s. for all $t \in [0,T]$.
  Applying $(I+\ep A)^{-1}$ to both sides we have, in view of Lemma
  \ref{l3},
  \[
  y^{(\ep)}(t) = e^{-tA}y^{(\ep)}_0 + \int_0^t e^{-(t-s)A} b^{(\ep)}(s)\,ds
  + \int_0^t e^{-(t-s)A} \Phi^{(\ep)}(s)\,dM(s),
  \]
  where we have used the notation $h^{(\ep)}:=(I+\ep A)^{-1}h$, for
  any $h$ for which it makes sense. In other words, $\eps{y}$ is the
  mild solution of
  \[
  d\eps{y}(t) + A\eps{y}(t)\,dt = \eps{b}(t)\,dt + \eps{\Phi}(t)\,dM(t),
  \qquad \eps{y}(0)=\eps{y}_0.
  \]
  On the other hand, notice that we can apply Lemma \ref{l2} since, by
  (\ref{eq:16}) and the Lipschitz continuity of $\Phi_0$,
  \[
  \E \int_0^T |\Phi_0(y(s))|^2_Q \,ds \lesssim 1 + \E \sup_{t \leq T}
  |y(t)|^2 < +\infty.
  \]
  Therefore, $\eps{y}$ is also a strong solution of the previous
  equation, for which we can apply It\^o's formula for the square of
  the norm, obtaining
  \begin{align*}
    \big|\eps{y}(t)\big|^2 &= \big|\eps{y}_0\big|^2
        + 2\int_0^t \bip{\eps{y}(s-)}{d\eps{y}(s)} + \big[\eps{y}\big](t)\\
    &= \big|\eps{y}_0\big|^2 - 2 \int_0^t \bip{A \eps{y}(s)}{\eps{y}(s)}\,ds
    + 2\int_0^t \bip{\eps{b}(s)}{\eps{y}(s)}\,ds\\
    &\quad + 2\int_0^t \ip{\eps{y}(s-)}{\eps{\Phi}(s)\,dM(s)} +
     \big[\eps{\Phi} \cdot M\big](t)\\
    &\leq |y_0|^2 + 2\int_0^t \bip{\eps{b}(s)}{\eps{y}(s)}\,ds
      + 2\int_0^t \ip{\eps{y}(s-)}{\eps{\Phi}(s)\,dM(s)}
      + [\Phi \cdot M](t),
  \end{align*}
  where the inequality follows by the monotonicity of $A$ and the
  contractivity of $(I+\varepsilon A)^{-1}$. Thus we also have
  \begin{align}
  \E\sup_{t \leq T} \big|\eps{y}(t)\big|^2 & \leq \E|y_0|^2 + 
  2 \E \sup_{t \leq T} \int_0^t \bip{\eps{b}(s)}{\eps{y}(s)}\,ds \nonumber \\
  & \qquad + 2 \E \sup_{t \leq T} \int_0^t
           \ip{\eps{y}(s-)}{\eps{\Phi}(s)\,dM(s)} + \E[\Phi \cdot M](T).
  \label{eq:56}
  \end{align}
  We are going to get a uniform (with respect to $\ep$) estimate for
  the third term on the right-hand side. For this, consider the local
  martingale
  \[
  N_\varepsilon(t) := \int_0^t \ip{\eps{y}(s-)}{\eps{\Phi}(s)\,dM(s)},
  \qquad t\in [0,T],
  \] 
  for which Davis' and Young's inequalities yield
  \begin{align*}
  \E \sup_{t \leq T} |N_\varepsilon(t)| &\lesssim  \E[N_\varepsilon](T)^{1/2}
  \leq  \E \sup_{s \leq T} |\eps{y} (s)| \,[\eps{\Phi} \cdot M](T)^{1/2}\\
  &\leq \frac{\alpha}{2} \, \E \sup_{s\leq T} |y(s)|^2
         + \frac{1}{2\alpha} \, \E[\Phi \cdot M](T),
  \end{align*}
  for any $\alpha>0$. On the other hand, similarly as above, we have
  \[
  \E[\Phi \cdot M] (T) = \E\langle \Phi \cdot M\rangle (T)
  = \E\int_0^T |\Phi_0(y(s))|^2_{Q_M} \,d\langle M\rangle(s)
  \lesssim 1 + \int_0^T \big(\E \sup_{r\leq s} |y(r)|^2 \big) ds,
  \]
  which yields
  \[
  \E \sup_{t \leq T} |N_\varepsilon(t)| \lesssim \frac{\al}{2} \,
  \E\sup_{s\leq T} |y(s)|^2 + \frac 1\al + \frac 1\al \int_0^T
  \big(\E \sup_{r\leq s} |y(r)|^2 \big) ds,
  \]
  hence also, by (\ref{eq:56}),
  \begin{align*}
  \E \sup_{t \leq T} \big|\eps{y}(t)\big|^2  & 
  \lesssim \E |y_0|^2 + \E \sup_{t\leq T} \int_0^t
    \bip{\eps{b}(s)}{\eps{y}(s)}\,ds\\
  & \qquad  + \frac 1\al + \frac{\al}{2} \, \E\sup_{s\leq T} |y(s)|^2
  + \frac 1\al \int_0^T \big(\E \sup_{r\leq s} |y(r)|^2\big) ds.
  \end{align*}
  Let us now pass to the limit as $\ep$ tends to zero in the previous
  inequality. We clearly have
  \[
  \lim_{\ep\to 0} \, \E\sup_{t\leq T} |\eps{y}(t)|^2 = 
  \E\sup_{t\leq T} |y(t)|^2,
  \]
  and, similarly,
  \[
  \lim_{\varepsilon \to 0}
  \E\sup_{t\leq T} \int_0^t \bip{\eps{b}(s)}{\eps{y}(s)}\,ds 
	= \E\sup_{t\leq T} \int_0^t \ip{f(y(s))}{y(s)} \,ds.
  \]
  In fact, one has $\bip{\eps{b}(s)}{\eps{y}(s)} \to
  \ip{f(y(s))}{y(s)}$ $\P$-a.s. as $\varepsilon \to 0$,
  and
  \[
  \big| \bip{\eps{b}(s)}{\eps{y}(s)} \big| 
  \leq |f(y(s))| \, |y(s)| 
  \lesssim 1+|y(s)|^2,
  \]
  with $\E\sup_{t\leq T} |y(t)|^2 < \infty$. We have thus proved
  \begin{align*}
  \E\sup_{t\leq T} \big|y(t)\big|^2 & \lesssim \E |y_0|^2  + \frac 1\al + \E\sup_{t\leq T} \int_0^t \ip{f(y(s))}{y(s)} \,ds \\
  & \qquad \quad    + \frac{\al}{2}
    \E\sup_{t\leq T} |y(t)|^2   + \frac 1\al \int_0^T \big(\E \sup_{r\leq s} |y(r)|^2 \big) ds,
  \end{align*}
  which implies (\ref{12}) choosing $\alpha$ small enough.
\end{proof}

We shall need the following integration-by-parts formula for
Hilbert-space-valued semimartingales (cf. \cite[\S 26.9]{Met}).
\begin{lemma}     \label{lem:ibpf}
  Let $Z_1, Z_2$ be two $X$-valued semimartingales. Then one has
  \begin{equation}     \label{eq:ibpf}
  \ip{Z_1(t)}{Z_2(t)} = \ip{Z_1(0)}{Z_2(0)} + \int_0^t \ip{Z_1(s-)}{dZ_2(s)}
  + \int_0^t \ip{Z_2(s-)}{dZ_1(s)} + [Z_1,Z_2](t),
  \end{equation}
  $\mathbb{P}$-a.s. for all $t\geq 0$, where $[Z_1,Z_2]$ denotes the
  quadratic (co)variation of $Z_1$ and $Z_2$.
\end{lemma}

Let us also recall, for the reader's convenience, a result of
A.~Jakubowski on weak convergence of stochastic integrals in infinite
dimensions (see \cite[Thm.~4]{Jakub-SSR96}).  Given a sequence
$\{Z^n\}_{n\in\enne}$ of $\mathcal{H}$-valued semimartingales on the
corresponding stochastic bases $(\Omega,\mathcal{F},\mathbf{F}^n,\P)$,
one says that $\{Z^n\}_{n\in\enne}$ satisfies the UT condition if, for
any sequence of $\mathbf{F}^n$-adapted $\mathcal{H}$-valued elementary
processes $\{\zeta^n\}_{n\in\enne}$ uniformly bounded by $1$, the
family of random variables
\[
 \Big\{ \int_0^T \ip{\zeta^n(s-)}{dZ^n(s)} \Big\}_{n\in\enne}
\]
is uniformly tight.
\begin{thm}     \label{thm:jaku}
For each $n \in \enne$, let $Z^n$ be a $\ch$-valued semimartingale
with respect to the stochastic basis
$(\Omega,\mathcal{F},\mathbf{F}^n,\P)$ and let $H^n$ be
$\mathbf{F}^n$-adapted and with paths in $D\big([0,T] \to
\mathcal{L}_2(\ch \to X)\big)$. If $\{Z^n\}_{n\in\enne}$ satisfies the
UT condition and $(H^n,Z^n) \to (H,Z)$ as $n \to \infty$ in
probability in $D\big([0,T] \to \mathcal{L}_2(\ch \to X) \times
\ch\big)$, then $Z$ is a semimartingale with respect to the natural
filtration generated by $(H,Z)$ and
\[
\big(H^n,Z^n,H^n_- \cdot Z^n\big) \to \big(H,Z,H_- \cdot Z\big)
\]
as $n \to \infty$ in probability in $D([0,T] \to \mathcal{L}_2(\ch \to
X) \times \ch \times X\big)$.
\end{thm}
Note that in \cite{Jakub-SSR96} the convergences are in law, not in
probability. However, everything goes through with convergence in
probability as well, just by inspection of the proof (as already
observed in a completely analogous setting in \cite[p.~1041]{KP-AP91},
cf. also \cite[Thm.~5.5]{KP2}).

Finally, a key role will be played by the following result of Br\'ezis
(see \cite[Thm.~18]{Bre-mm}):
\begin{thm}    \label{thm:brezis}
  Let $\beta$ be a maximal monotone graph in $\erre\times \erre$ such
  that $\dom(\beta)=\erre$ and $0\in \beta(0)$. Let
  $\{f_n\}_{n\in\enne}$ and $\{g_n\}_{n\in\enne}$ be sequences of real
  measurable functions defined on some finite measure space $(\Theta,
  \mathcal{A},\mu)$ such that $g_n \to g$ $\mu$-a.e. as $n \to
  \infty$, $f_n(x)\in \beta(g_n(x))$ for $\mu$-a.a. $x \in \Theta$,
  and $f_n \, g_n \in L^1(\Theta,\mu)$ with $\int_\Theta f_n \,
  g_n\,d\mu < C$ for all $n \in \enne$, with $C$ independent of
  $n$. Then there is a subsequence $\{n_k\}_{k\in\enne}$ such that
  $f_{n_k}$ converges to some $f$ in $L^1(\Theta,\mu)$ as $k \to
  \infty$, and $f(x) \in \beta(g(x))$ for $\mu$-a.a. $x \in \Theta$.
\end{thm}


\section{Proof of Theorem \ref{thm:main}} \label{proof:main}
\subsection{Some preparations}
Before starting to prove the main result, it is helpful to recall a
few well-known facts about mild solutions to stochastic wave equations
with Lipschitz nonlinearities. In particular, assuming just for the
purposes of this subsection that $\beta: \erre\to \erre$ is Lipschitz
continuous, (\ref{eq:sw}) can be written as the following evolution
equation on $H:=H^1_0 \times L^2$:
\begin{equation}     \label{eq:mw}
dU(t) + AU(t)\,dt + BU(t)\,dt = G(U(t-))\,dM(t),
\qquad U(0)=(u_0,v_0),
\end{equation}
where $U(t):=(u(t),v(t))$,
\begin{align*}
A: \dom(A) \subset H &\to H,\\
   (u,v) &\mapsto (-v,-\Delta u),
\end{align*}
$B: H \ni (u,v) \mapsto (0,\beta(u))$, and $G: H \ni (u,v) \mapsto
(0,G_0(u))$. The operator $-A$ generates a strongly continuous group
$S(t)=e^{-tA}$ on $H$, with
\[
S(t) = 
\begin{bmatrix}
\cos\big(t\sqrt{-\Delta}\big) &
\displaystyle \frac{1}{\sqrt{-\Delta}}\sin\big(t\sqrt{-\Delta}\big)\\
-\sqrt{-\Delta}\sin\big(t\sqrt{-\Delta}\big) &
\cos\big(t\sqrt{-\Delta}\big)
\end{bmatrix}
\]
(see e.g. \cite{Segal:wave}), and it is easily seen that (\ref{eq:mw})
admits a unique mild solution, i.e. an adapted c\`adl\`ag process $U$
such that the stochastic integral equation
\begin{equation}     \label{eq:mild}	
U(t) + \int_0^t S(t-s)B(U(s))\,ds = S(t) U(0)
+ \int_0^t S(t-s)G(U(s-))\,dM(s),
\end{equation}
is satisfied $\P$-a.s. for all $t \in [0,T]$. An elementary
computation based on the explicit form of $S(t)$ yields that
$U(t)=(u(t),v(t))$ satisfies (\ref{eq:mild}) if and only if it
satisfies the alternative form
\[
\left\{
\begin{aligned}
  &u(t) = u_0 + \int_0^t v(s)\,ds,\\
  &v(t) + \int_0^t S_{22}(t-s)\beta(u(s))\,ds = S_{21}(t)u_0 + S_{22}(t)v_0
  + \int_0^t S_{22}(t-s) G_0(u(s-))\,dM(s)
\end{aligned}
\right.
\]
$\P$-a.s. for all $t \in [0,T]$, where $S_{ij}(t)$ denotes the
$(i,j)$-th entry of the operator matrix $S(t)$. A further computation
(which is elementary, apart of having to appeal to a general
stochastic Fubini's theorem such as the one in \cite{Leon-Fub}) shows
that a mild solution satisfies the Duhamel's formulation
\[
u(t) + \int_0^t S_{12}(t-s) \beta(u(s)) \,ds = S_{11}(t) u_0 + S_{12}(t) v_0 
+ \int_0^t S_{12}(t-s) G_0(u(s-))\,dM(s).
\]
As already mentioned, this expression motivates the definition of weak
solution in the general case (i.e. without any Lipschitz assumption on
$\beta$).

\subsection{Proof of Theorem \ref{thm:main}}
Denoting the identity function by $I$, let
\[
\beta_\lambda = \frac{1}{\lambda}\big( I - (I+\lambda\beta)^{-1} \big),
\qquad \lambda>0,
\]
be the Yosida approximation of $\beta$. Recall that the maximal
monotonicity of $\beta$ implies that, for any $\lambda>0$,
$(I+\lambda\beta)^{-1}$ is a contraction defined on the whole real
line (see e.g. \cite[Prop.~2.2]{Bmax}). In particular, $\beta_\lambda$
is a function, not a graph. Furthermore, one can prove (see
e.g. \cite[Prop. 2.6]{Bmax}) that $\beta_\lambda:\erre\to\erre$ is
monotone, Lipschitz continuous with Lipschitz constant bounded above
by $2/\lambda$, and it satisfies $\beta_\lambda \in
\beta(I+\lambda\beta)^{-1}$.

Consider the regularized equation
\begin{equation}     \label{eq:6}
dU_\lambda(t)+ AU_\lambda(t)\,dt + B_\lambda(U_\lambda(t))\,dt
= G(U_\lambda(t-))\,dM(t),
\qquad U_\lambda(0)=(u_0,v_0),
\end{equation}
where $B_\lambda: (u,v) \mapsto (0,\beta_\lambda(u))$ is Lipschitz
continous from $H$ to itself. Then (\ref{eq:6}) admits a unique
c\`adl\`ag mild solution $U_\lambda$ such that
\[
 \E \sup_{t \leq T} |U_\lambda(t)|^2_H <+\infty
\]
(see e.g. \cite{Kote-Doob}). We are now going to establish a priori
estimates for $U_\lambda=(u_\lambda,v_\lambda)$.
\begin{prop}     \label{prop:ap}
  There exists a positive constant $C$, independent of $\lambda$, such
  that
  \[
  \E \sup_{t \leq T} \big( |u_\lambda(t)|_{H^1_0}^2 
  + |v_\lambda(t)|_{L^2}^2 \big) < C.
  \]
\end{prop}
\begin{proof}
  By Proposition \ref{prop3}, taking into account that $B_\lambda$ is
  Lipschitz continuous, we have
  \[
  \E \sup_{t \leq T} |U_\lambda(t)|^2 \lesssim 1 + \E|U_0|^2 + \E
  \sup_{t \leq T} \int_0^t
  \bip{-B_\lambda(U_\lambda(s))}{U_\lambda(s)} \,ds + \int_0^T \E
  \sup_{r\leq s} |U_\lambda(r)|^2 \, ds,
  \]
  where
  \[
  \int_0^t \bip{B_\lambda(U_\lambda(s))}{U_\lambda(s)}\,ds =
  \int_0^t \ip{\beta_\lambda(u_\lambda(s))}{v_\lambda(s)}_{L^2}\,ds.
  \]
  Let us introduce the Moreau-Yosida approximation of $j$ (recall that
  $j$ is a positive convex function such that $\beta=\partial j$),
  that is
  \[
  j_\lambda(x) := \inf_{y \in \erre}
  \big( j(y) + \frac{|x-y|^2}{2\lambda} \big), \qquad \lambda>0.
  \]
  Then one has (see e.g. \cite[Thm.~2.2.2]{barbu}) that $j_\lambda \in
  C^1(\erre)$, $\beta_\lambda=j'_\lambda$, and $j_\lambda \to j$
  pointwise as $\lambda \to 0$. Moreover, one also has $j_\lambda \leq
  j$ and, obviously, $j_\lambda \geq 0$. We can thus write, in view of
  the identity $u_\lambda(t) = u_0 + \int_0^t v_\lambda(s)\,ds$,
  recalling that $j'_\lambda=\beta_\lambda$,
  \begin{equation*}
    \int_0^t \ip{\beta_\lambda(u_\lambda(s))}{v_\lambda(s)}_{L^2}\,ds 
    = \int_D j_\lambda(u_\lambda(t,x))\,dx - \int_D
    j_\lambda(u(0,x))\,dx \geq -|j(u_0)|_{L^1}.
  \end{equation*}
  This implies
  \[
  \E \sup_{t \leq T} |U_\lambda(t)|^2 \lesssim 1 + \E|U_0|^2 
  + \E|j(u_0)|_{L^1}
  + \int_0^T \E \sup_{r\leq s} |U_\lambda(r)|^2 \,ds,
  \]
  which, by an application of Gronwall's inequality, yields the claim.
\end{proof}

\begin{prop}     \label{bound-beta}
  Let $\tilde{u}_\lambda(t):=(I+\lambda\beta)^{-1}u_\lambda(t)$, $t
  \in [0,T]$.  There exists a positive constant $C$, independent of
  $\lambda$, such that
  \[
  \E\int_0^T \ip{\beta(\tilde{u}_\lambda(s))}{\tilde{u}_\lambda(s)}_{L^2}
    \,ds < C.
  \]
\end{prop}
\begin{proof}
We split the proof in three steps. 
\smallskip\par\noindent
\emph{Step 1.} We introduce a regularized version of (\ref{eq:6})
admitting a strong solution. In particular, setting
\[
F_\lambda(t) := (0,f_\lambda(t)) := 
\big(0,\beta_\lambda(u_\lambda(t))\big),
\qquad 
\Gamma_\lambda(t) := (0,\gamma_\lambda(t)) := 
\big(0,G_0(u_\lambda(t-))\big)
\]
for all $t \in [0,T]$, it is clear that $U_\lambda$ is
the unique mild solution of
\begin{equation}     \label{eq:7}
dU_\lambda(t) + AU_\lambda(t)\,dt + F_\lambda(t)\,dt = \Gamma_\lambda(t)\,dM(t),
\qquad U_\lambda(0)=(u_0,v_0),
\end{equation}
or equivalently
\begin{equation}      \label{eq:9}
\left\{
\begin{aligned}
  &u_\lambda(t) = u_0 + \int_0^t v_\lambda (s)\,ds\\
  &v_\lambda(t) + \int_0^t S_{11}(t-s) f_\lambda(s)\,ds = S_{21}(t)u_0 + S_{11}v_0
  + \int_0^t S_{11}(t-s) \gamma_\lambda(s)\,dM(s).
\end{aligned}
\right.
\end{equation}
Using the notation $\eps{h}:=(I-\varepsilon\Delta)^{-1}h$ for any
``object'' $h$ for which the expression makes sense, we may write,
recalling that $(I+\varepsilon\Delta)^{-1}$ and $S_{ij}(\cdot)$,
$i,\,j=1,\,2$, commute,
\begin{equation}     \label{eq:xx}
\left\{
\begin{aligned}
&\eps{u}_\lambda(t) = \eps{u_0} + \int_0^t \eps{v}_\lambda(s)\,ds,\\
&\eps{v}_\lambda(t) + \int_0^t S_{11}(t-s) \eps{f}_\lambda(s)\,ds 
= S_{21}(t) \eps{u}_0 + S_{11} \eps{v}_0 
  + \int_0^t S_{11}(t-s) \eps{\gamma}_\lambda(s)\,dM(s),
\end{aligned}
\right.
\end{equation}
or equivalently, $\eps{U}_\lambda(t)=(\eps{u}_\lambda(t),\eps{v}_\lambda(t))$
is the unique mild solution of
\begin{equation}     \label{eq:10}
d\eps{U}_\lambda(t)+ A\eps{U}_\lambda(t)\,dt +
\eps{F}_\lambda(t)\,dt = \eps{\Gamma}_\lambda(t)\,dM(t),
\qquad 
\eps{U}_\lambda(0)=(\eps{u}_0,\eps{v}_0),
\end{equation}
where
\[
\eps{F}_\lambda(t) = \big(0,\eps{f}_\lambda(t)\big),
\qquad
\eps{\Gamma}_\lambda(t) = \big(0,\eps{\gamma}_\lambda(t)\big)
\]
for all $t\in [0,T]$. By Lemma \ref{l2} we actually have that
$\eps{U}_\lambda$ is a strong solution of (\ref{eq:10}).
\smallskip\par\noindent
\emph{Step 2.} 
Let $\eps{V_\lambda}(t):=(0,(I+\lambda\beta)^{-1} \eps{u}_\lambda(t))$, $t \in
[0,T]$. Then both $\eps{U}_\lambda$, as a strong solution of
(\ref{eq:10}), and $\eps{V}_\lambda$, are $H$-valued semimartingales, for
which the integration-by-parts formula (\ref{eq:ibpf}) yields
\begin{align*}
\bip{\eps{U}_\lambda(t)}{\eps{V}_\lambda(t)} & =
\bip{\eps{U}_\lambda(0)}{\eps{V}_\lambda(0)}
+ \int_0^t \bip{\eps{U}_\lambda(s-)}{d\eps{V}_\lambda(s)}
\\
&\qquad + \int_0^t \bip{\eps{V}_\lambda(s-)}{d\eps{U}_\lambda(s)} + \big[\eps{U}_\lambda,\eps{V}_\lambda\big](t).
\end{align*}
Taking into account (\ref{eq:10}) and the definitions of $\eps{U}_\lambda$,
$\eps{V}_\lambda$, $\eps{F}_\lambda$, $\eps{\Gamma}_\lambda$, we obtain
\begin{align}
&\E\bip{\eps{v_\lambda}(t)}{(I+\lambda\beta)^{-1}(\eps{u_\lambda}(t))}_{L^2} = 
\E\bip{\eps{v}_0}{(I+\lambda\beta)^{-1} \eps{u}_0}_{L^2} \nonumber\\
&\qquad + \E\int_0^t \bip{\eps{v_\lambda}(s)}
        {d\big( (I+\lambda\beta)^{-1}(\eps{u_\lambda}(s))\big)}_{L^2}\nonumber\\
&\qquad + \E\int_0^t \bip{(I+\lambda\beta)^{-1} \eps{u_\lambda}(s)}
                   {\Delta \eps{u_\lambda}(s)}_{L^2} \,ds \nonumber\\
&\qquad - \E\int_0^t \bip{(I+\lambda\beta)^{-1} \eps{u_\lambda}(s)}
           {(I-\ep\Delta)^{-1}\beta_\lambda(u_\lambda(s))}_{L^2} \,ds \nonumber\\
&\qquad + \E\int_0^t \bip{(I+\lambda\beta)^{-1}\eps{u_\lambda}(s-)}
                    {\eps{\gamma}_\lambda(s)\,dM(s)}_{L^2} \,ds \nonumber\\
&\qquad + \E\big[\eps{v}_\lambda,(I+\lambda\beta)^{-1} \eps{u}_\lambda\big](t).
\label{eq:11}
\end{align}
Note that, by (\ref{eq:xx}),
$\partial_s\eps{u}_\lambda(s)=\eps{v}_\lambda(s)$, and
\begin{align*}
(I-\lambda\beta)^{-1} \eps{u}_\lambda(s) &= \eps{u}_\lambda(s)
+ (I-\lambda\beta)^{-1}\big(\eps{u}_\lambda(s)\big) - \eps{u}_\lambda(s)\\
&= \eps{u}_\lambda(s) - \lambda\beta_\lambda\big(\eps{u}_\lambda(s)\big),\\
\end{align*}
hence
\begin{align*}
&\E\int_0^t \bip{\eps{v}_\lambda(s)}
            {d\big( (I+\lambda\beta)^{-1}(\eps{u}_\lambda(s)) \big)}_{L^2} \, ds \\
&\qquad = \E\int_0^t \bip{\eps{v}_\lambda(s)}{d\big(
          \eps{u}_\lambda(s) - \lambda\beta_\lambda(\eps{u}_\lambda(s))
          \big)}_{L^2}\\
&\qquad = \E\int_0^t |\eps{v}_\lambda(s)|^2_{L^2} \,ds
- \lambda \E\int_0^t \bip{\eps{v}_\lambda(s)}
{d\beta_\lambda(\eps{u}_\lambda(s))}_{L^2}.
\end{align*}
Denoting a family of $C^\infty_c$ mollifiers by
$\{\zeta_\delta\}_{\delta>0}$, and setting
$\beta_{\lambda\delta}:=\beta_\lambda\ast\zeta_\delta$, one has
\[
\E\int_0^t \bip{\eps{v}_\lambda(s)}
{d\beta_{\lambda\delta}(\eps{u}_\lambda(s))}_{L^2}
\xrightarrow{\delta\to 0}
\E\int_0^t \bip{\eps{v}_\lambda(s)}
{d\beta_\lambda(\eps{u}_\lambda(s))}_{L^2}
\]
and
\[
\E\int_0^t \bip{\eps{v}_\lambda(s)}
{d\beta_{\lambda\delta}(\eps{u}_\lambda(s))}_{L^2}
= \E\int_0^t \bip{\eps{v}_\lambda(s)}
{\beta'_{\lambda\delta}(\eps{u}_\lambda(s))\eps{v}_\lambda(s)}_{L^2}\,ds
\geq 0,
\]
because, recalling that $\beta_\lambda$ is a monotonically increasing
function, $\beta'_{\lambda\delta}(x) = \beta'_\lambda \ast
\zeta_\delta(x) \geq 0$ for all $x \in \erre$. Therefore we have
\[
\E\int_0^t \bip{\eps{v}_\lambda(s)}
{d\big( (I+\lambda\beta)^{-1}(\eps{u}_\lambda(s)) \big)}_{L^2} \, ds
\leq \E\int_0^t |\eps{v}_\lambda(s)|^2_{L^2} \,ds.
\]

Similarly, we have
\begin{align*}
& \E\int_0^t \bip{(I+\lambda\beta)^{-1}(\eps{u_\lambda}(s))}
{\Delta \eps{u_\lambda}(s)}_{L^2}\,ds\\
&\qquad = \E\int_0^t\bip{\eps{u_\lambda}(s)}{\Delta\eps{u_\lambda}(s)}_{L^2}\,ds
- \E\int_0^t \lambda\bip{\beta_\lambda(\eps{u_\lambda}(s))}
                     {\Delta\eps{u_\lambda}(s)}_{L^2}\,ds\\
&\qquad = -\E\int_0^t \bip{\nabla\eps{u_\lambda}(s)}
                            {\nabla\eps{u_\lambda}(s)}_{L^2}\,ds
+ \E\int_0^t \lambda \bip{\beta_\lambda'(\eps{u_\lambda}(s)) \nabla\eps{u_\lambda}(s)}
             {\nabla\eps{u_\lambda}(s)}_{L^2}\,ds\\
&\qquad \leq 2\E\int_0^t |\nabla\eps{u_\lambda}(s)|^2_{L^2}\,ds,
\end{align*}
where the term involving $\beta'_\lambda$ can be interpreted, as
above, as limits of more regular expressions obtained replacing
$\beta_\lambda$ with $\beta_{\lambda\delta}$.

Moreover, by an argument completely similar to one used in the proof
of Proposition \ref{prop3}, we get
\[
\E \sup_{t \leq T}  \int_0^t \bip{(I+\lambda\beta)^{-1}\eps{u_\lambda}(s-)}
                    {\eps{\gamma}_\lambda(s)\,dM(s)}_{L^2} \,ds < +\infty,
\]
which implies that the stochastic integral appearing on the right-hand side
of (\ref{eq:11}) is a martingale, hence with expectation zero. 
Finally, we have that
\[
\big[\eps{v}_\lambda,(I+\lambda\beta)^{-1} \eps{u}_\lambda\big](t)=0
\]
for all $t \in [0,T]$ because, as it follows by (\ref{eq:xx}),
\[
(I+\lambda\beta)^{-1} \eps{u_\lambda}(t)= (I+\lambda\beta)^{-1} \eps{u}_0
+ \int_0^t (I+\lambda\beta)^{-1} \eps{v_\lambda}(s) \, ds
\]
is a process with finite variation and continuous paths. We have thus
proved that
\begin{align*}
&\E\int_0^T \bip{(I+\lambda\beta)^{-1}\eps{u}_\lambda(s)}
   {(I-\ep\Delta)^{-1}\beta_\lambda(u_\lambda(s))}_{L^2} \,ds\\
& \qquad \leq \E \bip{\eps{v}_0}{(I+\lambda\beta)^{-1}\eps{u}_0}_{L^2}
+ \E \big|\bip{\eps{v_\lambda}(T)}{(I+\lambda\beta)^{-1}\eps{u_\lambda}(T)}_{L^2}
     \big|\\
& \qquad\quad + \E\int_0^T \big|\eps{v}_\lambda(s)\big|^2_{L^2}\,ds
+ 2\E\int_0^T \big|\nabla\eps{u}_\lambda(s)\big|^2_{L^2}\,ds.
\end{align*}
Furthermore, Cauchy-Schwarz' inequality and the contractivity of
$(I+\lambda\beta)^{-1}$ and $(I-\varepsilon\Delta)^{-1}$ yield
\begin{align*}
\E \bip{\eps{v}_0}{(I+\lambda\beta)^{-1}\eps{u}_0}_{L^2} &\leq
\big(\E|v_0|^2_{L^2}\big)^{1/2} \, \big(\E|u_0|^2_{L^2}\big)^{1/2},\\
\E\big|\bip{\eps{v_\lambda}(T)}{(I+\lambda\beta)^{-1}\eps{u_\lambda}(T)}_{L^2}\big|
&\leq \big(\E|v_\lambda(T)|^2_{L^2}\big)^{1/2} \,
     \big(\E|u_\lambda(T)|^2_{L^2}\big)^{1/2}\\
&\leq \left( \E \sup_{t \leq T} |v_\lambda(t)|^2_{L^2} \right)^{1/2} \,
      \left( \E \sup_{t \leq T} |u_\lambda(t)|^2_{L^2} \right)^{1/2},\\
\E\int_0^T \big|\eps{v}_\lambda(s)\big|^2_{L^2}\,ds
&\leq \E\int_0^T \big|v_\lambda(s)\big|^2_{L^2}\,ds
\leq T \, \E \sup_{t\leq T} |v_\lambda(t)|^2_{L^2}, \\
\E\int_0^T \big|\nabla\eps{u}_\lambda(s)\big|_{L^2}\,ds
& \leq \E\int_0^T \big|\nabla u_\lambda(s)\big|_{L^2}\,ds
\leq T \, \E \sup_{t \leq T} |u_\lambda(t)|^2_{H^1_0}.
\end{align*}
Appealing to Proposition \ref{prop:ap}, we infer that there exists a
constant $C$, independent of $\varepsilon$ and of $\lambda$, such that
\begin{equation}     \label{eq:bacca}
\E\int_0^t \bip{(I+\lambda\beta)^{-1}\eps{u}_\lambda(s)}
   {(I-\ep\Delta)^{-1}\beta_\lambda(u_\lambda(s))}_{L^2} \,ds
< C.
\end{equation}
\smallskip\par\noindent
\emph{Step 3.} We shall now pass to the limit as $\varepsilon \to 0$
in the last inequality.  Since the operator $(I+\lambda\beta)^{-1}$ is
bounded, we have that
\[
\lim_{\ep\rightarrow 0}
\bip{(I+\lambda\beta)^{-1} \eps{u_\lambda}(s)}
{(I-\ep\Delta)^{-1} \beta_\lambda(u_\lambda(s))}_{L^2} 
= \bip{(I+\lambda\beta)^{-1} u_\lambda(s)}{\beta_\lambda(u_\lambda(s))}_{L^2}
\]
$\P$-a.s. for all $s \in [0,T]$. On the other hand, the contractivity
of $(I+\lambda\beta)^{-1}$ and $(I-\ep\Delta)^{-1}$, the Lipschitz
continuity of $\beta_\lambda$, and Cauchy-Schwarz' inequality yield
\begin{align*}
&\E\int_0^T \big| \big\langle (I+\lambda\beta)^{-1} \eps{u_\lambda}(s),
(I-\ep\Delta)^{-1} \beta_\lambda(u_\lambda(s)) \big\rangle_{L^2}\big| \,ds \\
&\qquad\qquad \lesssim \E\int_0^T (1+|u_\lambda(s)|^2_{L^2})\,ds < \infty.
\end{align*}
Hence, by the dominated convergence theorem, we obtain
\begin{align*}
& \E\int_0^T \big\langle (I+\lambda\beta)^{-1}
    u_\lambda(s), \beta_\lambda(u_\lambda(s)) \big\rangle_{L^2} \,ds\\
&\qquad\qquad = \lim_{\ep \to 0} \E\int_0^T 
\big\langle(I+\lambda\beta)^{-1} \eps{u_\lambda}(s),
(I-\ep\Delta)^{-1} \beta_\lambda(u_\lambda(s)) \big\rangle_{L^2} \,ds
< C,
\end{align*}
where $C$ is the same constant (independent of $\lambda$) appearing in
(\ref{eq:bacca}), thus completing the proof, upon recalling that
$\beta_\lambda \in \beta(I+\lambda\beta)^{-1}$.
\end{proof}

\begin{proof}[Proof of Theorem \ref{thm:main}]
Let $(u_\lambda,v_\lambda)$ be the
solution of the regularized equation (\ref{eq:6}), for which one has
\begin{equation}     \label{eq:urca}
u_\lambda(t) = u_0 + \int_0^t v_\lambda(s)\,ds
\end{equation}
for all $t\in [0,T]$, and, by Proposition \ref{prop:ap},
\[
\E \sup_{t \leq T} |v_\lambda(t)|^2_{L^2} < C
\]
for a constant $C$ that does not depend on $\lambda$.
This implies, for any $\eta>0$,
\[
\P\big( \sup_{t\leq T} |v_\lambda(t)|^2_{L^2} > \eta \big)
\leq \frac{1}{\eta} \E \big( \sup_{t\leq T} |v_\lambda(t)|^2_{L^2} \big)
< \frac{C}{\eta},
\]
thus also
\begin{equation}     \label{eq:rebo2}
\lim_{\eta \to \infty} \sup_{\lambda}
\P\big( \sup_{t\leq T} |v_\lambda(t)|^2_{L^2} > \eta \big) = 0.
\end{equation}
Let us now show that, for any $t \in [0,T]$ and $\delta>0$,
there exists a compact subset $K=K(t,\delta)$ of $L^2$ such that
$\P(u_\lambda(t) \in K) > 1-\delta$. In fact, denoting by $B_R$ the
ball of radius $R$ centered at the origin of $H^1_0$ and recalling
that $H^1_0$ is compactly embedded in $L^2$, we have that $B_R$ is a
compact subset of $L^2$, and
\begin{align*}
\P\big( u_\lambda(t) \in B_R \big) =
1 - \P\big( |u_\lambda(t)|_{H^1_0} > R \big)
\geq 1 - \frac{1}{R^2} \E\sup_{t\leq T} |u_\lambda(t)|^2_{H^1_0}
> 1 - \frac{C}{R^2},
\end{align*}
i.e. it is enough to choose $R=\sqrt{C/\delta}$ and $K=B_R$. This
observation and (\ref{eq:rebo2}) imply, thanks to a corollary to a
theorem of Rebolledo (see e.g. \cite[\S II.4.4]{Met-SNS}), that
$\{u_\lambda\}_\lambda$ is uniformly tight in $D([0,T] \to L^2)$.
By Skorohod's representation theorem (see
e.g. \cite[\S~8.5]{Bog-MT2}), there exists a probability space
$(\bar{\Omega},\bar{\mathcal{F}},\bar{\P})$ and a sequence of random
vectors $\xi_n:=(\bar{u}_n,\bar{M}_n) \in D([0,T] \to L^2\times L^2)$
such that $\xi_n \to \xi:=(\bar{u},\bar{M})$ $\bar{\P}$-a.s. as $n \to
\infty$, and the laws of $\xi_n$ and $(u_{\lambda_n},M)$ coincide for
each $n$, for some subsequence $\{\lambda_n\}_{n\in\enne}$ of
$\lambda$. Let us set
\[
\mathcal{B}_n = \big( \bar{\Omega},\bar{\cf},\bar{\mathbf{F}}_n,\bar{\P} \big),
\qquad
\mathcal{B} = \big( \bar{\Omega},\bar{\cf},\bar{\mathbf{F}},\bar{\P} \big),
\]
where $\bar{\mathbf{F}}_n$ and $\bar{\mathbf{F}}$ are the filtrations
generated by
$(\bar{u}_n,\bar{M}_n)$ and $(\bar{u},\bar{M})$, respectively. Then
$\bar{u}_n$ is $\bar{\mathbf{F}}_n$-adapted and c\`adl\`ag, and $\bar{u}$ is
$\bar{\mathbf{F}}$-adapted and c\`adl\`ag, since $\bar{u}_n$ converges
$\bar{\P}$-a.s. to $\bar{u}$ as $n \to \infty$ in the Skorohod topology.

Let us assume, for the time being, that the process $\bar{M}_n$ is a
$\bar{\mathbf{F}}_n$-martingale for each $n$. Then the process
$\bar{M}$ is a $\bar{\mathbf{F}}$-martingale by a slight modification
of the proof of \cite[Prop.~IX.1.10]{JacShi}, taking into account
[op.~cit., Prop.~IX.1.12 and Rmk.~VI.1.10], as well as the obvious
inequality
\[
\sup_{t\leq T} \bar{\E} |\bar{M}_n(t)|^2 = 
\sup_{t\leq T} \E|M(t)|^2 < \infty
\]
(where $\bar{\E}$ denotes expectation with respect to $\bar{\P}$ on the
stochastic basis $\mathcal{B}_n$), which implies that $\bar{M}_n$ is
uniformly integrable. On the other hand, a completely similar argument
proves that $\bar{M}_n$ is indeed a $\mathbf{F}_n$-martingale for each
$n$: fix $n$ and set $M^k:=M$ and $Y^k:=u_n$, $k \in \enne$, so that
$(Y^k,M^k)$ trivially converges in law to $(\bar{u}_n,\bar{M}_n)$,
which has the same law of $(u_n,M)$. Since $\{M^k\}_{k\in\enne}$ is
obviously uniformly integrable, we conclude that $\bar{M}_n$ is a
martingale with respect to the filtration generated by
$(\bar{u}_n,\bar{M}_n)$, as it follows by the above mentioned results
of \cite{JacShi}.

Setting $\beta_n:=\beta_{\lambda_n}$, we also have that, for each
fixed $n$, $\bar{u}_n$ solves $\bar{\P}$-a.s. the integral equation
\begin{multline}     \label{eq:57}
\bar{u}_n(t)
+ \int_0^t \frac{1}{\sqrt{-\Delta}} \sin((t-s)\sqrt{-\Delta})
           \beta_n(\bar{u}_n(s))\,ds \\
= \cos(t\sqrt{-\Delta}) u_0
+ \frac{1}{\sqrt{-\Delta}} \sin(t\sqrt{-\Delta}) v_0
+ \int_0^t \frac{1}{\sqrt{-\Delta}} \sin((t-s)\sqrt{-\Delta})
           G_0(\bar{u}_n(s-))\,d\bar{M}_n(s),
\end{multline}
because the same equation is satisfied with $\bar{u}_n$ and
$\bar{M}_n$ replaced by $u_{\lambda_n}$ and $M$, respectively (for
instance by an argument such as the one used in \cite[p.~89]{krylov}).

\begin{lemma}     \label{l1}
  There is a subsequence $\{n_k\}_{k\in \enne}$ such that
  $\bar{u}_{n_k}\to \bar{u}$ $\bar{\P} \otimes \mathrm{Leb}$-a.e. on
  $\bar{\Omega} \times [0,T] \times D$.
\end{lemma}
\begin{proof}
  Note that we clearly have
  \begin{equation}    \label{eq:63} 
    \bar{\E} \int_0^T |\bar{u}_n(t)-\bar{u}(t)|_{L^1} \, dt \lesssim
    \bar{\E} \int_0^T |\bar{u}_n(t)-\bar{u}(t)|_{L^2} \, dt.
  \end{equation}
  Since (\ref{eq:urca}) implies $\bar{u}_n \in C([0,T]\to L^2)$, upon
  recalling that the Skorokhod topology on $D([0,T] \to L^2)$ induces
  the uniform topology on its subspace $C([0,T]\to L^2)$, we have that
  $\bar{u}_n$ converges to $\bar{u}$ $\bar{\P}$-a.s. in $C([0,T]\to
  L^2)$. That is
  \begin{equation}
    \sup_{t\leq T} |\bar{u}_n(t)-\bar{u}(t)|_{L^2} \to 0
    \label{eq:64}
  \end{equation}
  $\bar{\P}$-a.s. as $n \to \infty$. Set $X_n:= \int_0^T
  |\bar{u}_n(t)-\bar{u}(t)|_{L^2} \, dt$. By (\ref{eq:64}), the
  sequence $X_n$ converges to $0$ in probability. Assume for the
  moment that
  \begin{equation}
    \sup_n \bar{\E} |X_n|^2_{L^2} <+\infty.
    \label{eq:65}
  \end{equation}
  Then $\{X_n\}_{n\in\enne}$ is uniformly integrable and $\bar{\E}X_n
  \to 0$ as $n \to \infty$, whence, by (\ref{eq:63}), $\lim_{n \to
    \infty} \bar{u}_n = \bar{u}$ in $L^1(\bar{\Omega}\times D \times
  [0,T])$, thus also $\bar{u}_{n_k} \to \bar{u}$ $\bar{\P} \otimes
  \mathrm{Leb}$-a.e. on $\bar{\Omega} \times D \times [0,T]$ as $k \to
  \infty$, along a subsequence $\{n_k\}_{k\in\enne}$.

  It remains to justify (\ref{eq:65}). Observe that
  \begin{align*}
  \bar{\E} |X_n|^2 & \lesssim \bar{\E} \int_0^T
                     |\bar{u}_n(t)-\bar{u}(t)|^2_{L^2} \, dt\\
  & \lesssim \bar{\E} \int_0^T |\bar{u}_n(t)|^2_{L^2}
    + \bar{\E} \int_0^T |\bar{u}(t)|^2_{L^2} \, dt\\
  & =: I_1 + I_2,
  \end{align*}
  and
  \[
  I_1 \lesssim \sup_n \sup_{t\leq T} \bar{\E} |\bar{u}_n(t)|^2_{L^2} < +\infty
  \]
  by Proposition \ref{prop:ap} and the fact that $\bar{u}_n$ has the
  same law as $u_{\lambda_n}$. On the other hand, (\ref{eq:64}) implies
  \[
  \int_0^T |\bar{u}_n(t)|^2_{L^2}\, dt \to 
  \int_0^T |\bar{u}(t)|^2_{L^2}\, dt
  \]
  $\bar{\P}$-a.s. as $n \to \infty$, hence, by Fatou's lemma,
  \[
  I_2 = \bar{\E} \int_0^T |\bar{u}(t)|^2_{L^2} \,dt \leq 
  \liminf_n \bar{\E} \int_0^T |\bar{u}_n(t)|^2_{L^2} \,dt
  \lesssim \sup_n \sup_{t\leq T} \bar{\E} |\bar{u}_n(t)|^2_{L^2}
  < +\infty,
  \]
  which concludes the proof.
\end{proof}
Let us now consider the convergence of the stochastic integrals in
(\ref{eq:57}).
\begin{lemma}     \label{stochint}
  One has, for each $t \in [0,T]$,
  \[
  \int_0^t S_{12}(t-s) G_0(\bar{u}_n(s-))\,d\bar{M}_n(s) \to
  \int_0^t S_{12}(t-s) G_0(\bar{u}(s-))\,d\bar{M}(s)
  \]
  in probability as $n \to \infty$.
\end{lemma}
\begin{proof}
  By a basic trigonometric identity we can write
  \[
  \sin((t-s)\sqrt{-\Delta}) = \sin(t \sqrt{-\Delta})
  \cos(-s \sqrt{-\Delta}) + \cos(t \sqrt{-\Delta}) \sin(-s \sqrt{-\Delta}).
  \]
  Therefore, setting
  \[
  H_n^1(s) := \cos(-s \sqrt{-\Delta}) G_0(\bar{u}_n(s)), \qquad
  H_n^2(s) := \sin(-s \sqrt{-\Delta}) G_0(\bar{u}_n(s)),
  \]
  for all $s \in [0,T]$,
  we get
  \begin{equation}     \label{eq:43}
  \begin{split}
  \int_0^t S_{12}(t-s) G_0(\bar{u}_n(s-)) \,d\bar{M}_n(s) &=
  S_{12}(t) \int_0^t H^1_n(s-)\,d\bar{M}_n(s)\\
  &\quad + \frac{1}{\sqrt{-\Delta}} \cos( t\sqrt{-\Delta})
  \int_0^t H^2_n(s-)\,d\bar{M}_n(s).
  \end{split}
  \end{equation}
  By the continuity of $G_0$ and the boundedness of
  $\cos(-s\sqrt{-\Delta})$ and $\sin(-s \sqrt{-\Delta})$, we infer
  $(H^1_n,\bar{M}_n) \to (H^1,\bar{M})$ and $(H^2_n,\bar{M}_n) \to
  (H^2,\bar{M})$ in $D([0,T] \to \mathcal{L}^Q_2 \times
  L^2)$ in probability as $n \to \infty$, where
  \[
  H^1(s) := \cos(-s\sqrt{-\Delta}) G_0(\bar{u}(s)),
  \qquad 
  H^2(s) := \sin(-s\sqrt{-\Delta}) G_0(\bar{u}(s))
  \]
  for all $s \in [0,T]$. Let us now show that the sequence
  $\{\bar{M}_n\}_{n\in\enne}$ satisfies the UT condition specified
  just before the statement of Theorem \ref{thm:jaku}.  In fact,
  denoting a sequence of elementary processes as in \S\ref{sec:aux},
  mutatis mutandis, by $\{\zeta_n\}_{n\in\enne}$, (\ref{eq:16})
  and Markov's inequality yield
  \[
    \bar{\P}\left( \left| \int_0^t \langle \zeta_n(s-),d
        \bar{M}_n(s)\rangle_{\ch}\right| > \eta \right)
    \leq  T \, \frac{1}{\eta^2}\, \bar{\E} \sup_{s\leq t} |\zeta_n(s)|^2_{\ch}
    \leq \frac{T}{\eta^2},
  \]
  which implies, observing that the upper bound just obtained does not
  depend on $n$,
  \[
  \lim_{\eta\to \infty} \sup_n \bar{\P}\left( \left| \int_0^t \langle
      \zeta_n(s-), d \bar{M}_n(s)\rangle_{\ch}\right| > \eta \right) = 0
  \]
  for all $t \leq T$, i.e. the UT condition is verified.
  At this point we can apply Jakubowski's result Theorem \ref{thm:jaku}
  to deduce that
  \begin{align*}
  \int_0^\cdot H^n_1(s-) \,d\bar{M}_n(s) &\to \int_0^\cdot \cos(-s
  \sqrt{-\Delta}) G_0(\bar{u}(s-)) \,d\bar{M}(s),\\
  \int_0^\cdot H^n_2(s-) \,d\bar{M}_n(s) &\to
  \int_0^\cdot \sin(-s \sqrt{-\Delta}) G_0(\bar{u}(s-)) \,d\bar{M}(s),
  \end{align*}
  in probability in $D([0,T] \to L^2)$ as $n \to \infty$, which allows
  to conclude, in view of (\ref{eq:43}).
\end{proof}
We can now conclude the proof of Theorem \ref{thm:main}.  By
Proposition \ref{bound-beta} and Skorohod's representation, there
exists a constant $C$, independent of $n$, such that
\[
\bar{\E} \int_0^T \ip{\beta(J_n \bar{u}_n)}{J_n \bar{u}_n}_{L^2} \,ds < C,
\]
where $J_n := (I+\lambda_n\beta)^{-1}$. Therefore, by Theorem
\ref{thm:brezis}, we have that there exists $g \in \beta(\bar{u})$
such that, on a further subsequence, still denoted by $n$,
\[
\beta_n(\bar{u}_n) \to g
\qquad \text{in }  
L^1\big(\bar{\Omega}\times[0,T]\times D,\bar{\P}\otimes\mathrm{Leb}\big)
\]
as $n \to \infty$. In particular, passing to a further subsequence if
necessary,
\begin{equation}     \label{eq:44}
\beta_n(\bar{u}_n) \to g \qquad \text{in }  L^1([0,T] \to L^1)
\end{equation}
$\bar{\P}$-a.s. as $n \to \infty$.
Let us define the scale of Hilbert spaces
\[
\mathfrak{H}^m := \mathfrak{H}^m_1 \times \mathfrak{H}^m_2 := 
\dom\big((I-\Delta)^{m/2}\big) \times
\dom\big((I-\Delta)^{(m-1)/2}\big),
\qquad m \in \erre,
\]
where, for each $m \in \erre$, $\dom\big((I-\Delta)^{m/2}\big)$ is
endowed with the norm
\[
\| x \|_m = \big| (I-\Delta)^{m/2} x \big|_{L^2}.
\]
For instance, $\mathfrak{H}^0=L^2 \times H^{-1}$ and
$\mathfrak{H}^1=H^1_0 \times L^2$, which are the traditional Hilbert
spaces on which the strongly continuous group $\{S(t)\}_{t\in\erre}$
associated to the linear wave equation is considered. Since, for each
$s \in \erre$, $(I-\Delta)^{s/2}: \mathfrak{H}^m \to
\mathfrak{H}^{m-s}$ is an isometric isomorphism (considered
componentwise) and $(I-\Delta)^{s/2}$ commutes with $S(t)$ for any
$s$, $t$, one immediately verifies that $\{S(t)\}_{t\in\erre}$ can be
extended (or restricted) to a strongly continuous group on
$\mathfrak{H}^m$ for all $m \in \erre$ (cf. \cite[\S II.5]{EnNa} for a
related general scheme to extend semigroups of operators to so-called
Sobolev towers).  By classical Sobolev embedding theorems (see
e.g. \cite{Mazya}), there exists $m>0$ such that $L^1 \subset
\dom\big((I-\Delta)^{-(m+1)/2}\big)$ with continuous embedding,
therefore, in view of (\ref{eq:44}),
\[
\beta_n(\bar{u}_n) \to g \qquad \text{in }
L^1\big( [0,T] \to \mathfrak{H}^{-m}_2 \big)
\]
$\bar{\P}$-a.s. as $n \to \infty$, and so
\[
\int_0^t S_{12}(t-s) \beta_n(\bar{u}_n(s))\,ds 
\to \int_0^t S_{12}(t-s)g(s) \,ds
\]
$\bar{\P}$-a.s. as $n \to \infty$ by continuity of
$\{S_{12}(t)\}_{t\in\erre}$ in $\mathfrak{H}^{-m}_2$.

Summing up, we have obtained that $(\bar{u},g)$ solves
$\bar{\P}$-a.s. the equation
\begin{align*}
& \bar{u}(t) + \int_0^t \frac{1}{\sqrt{-\Delta}} \sin((t-s)\sqrt{-\Delta}) g(s)\,ds \\
& = \cos(t\sqrt{-\Delta}) u_0
+ \frac{1}{\sqrt{-\Delta}} \sin(t\sqrt{-\Delta}) v_0 
+ \int_0^t \frac{1}{\sqrt{-\Delta}} \sin((t-s)\sqrt{-\Delta})
           G_0(\bar{u}(s-))\,d\bar{M}(s),
\end{align*}
where all random vectors have to be considered as taking values in
$\mathfrak{H}^{-m}_2$, thus proving that $\bar{u}$ is a
(probabilistically) weak solution of (\ref{eq:fw}), in the sense of
Definition \ref{def:wsol}.
\end{proof}

\bibliographystyle{amsplain}

\providecommand{\bysame}{\leavevmode\hbox to3em{\hrulefill}\thinspace}
\providecommand{\MR}{\relax\ifhmode\unskip\space\fi MR }
\providecommand{\MRhref}[2]{%
  \href{http://www.ams.org/mathscinet-getitem?mr=#1}{#2}
}
\providecommand{\href}[2]{#2}

\end{document}